\documentclass[12pt]{article}%
\usepackage{amsfonts}
\usepackage{sw20bams}
\usepackage{amsmath}
\usepackage{amssymb}
\usepackage{graphicx}%
\providecommand{\U}[1]{\protect\rule{.1in}{.1in}}
\begin{document}

\title{Generalized Logistic Maps and Convergence}
\author{Steven Finch}
\date{July 15, 2025}
\maketitle

\begin{abstract}
We treat three cubic recurrences, two of which generalize the famous iterated
map $x\mapsto x(1-x)$ from discrete chaos theory. \ A\ feature of each
asymptotic series developed here is a constant, dependent on the initial
condition but otherwise intrinsic to the function at hand.

\end{abstract}

\footnotetext{Copyright \copyright \ 2025 by Steven R. Finch. All rights
reserved.}The quadratic recurrence%
\[%
\begin{tabular}
[c]{llll}%
$x_{0}=\dfrac{1}{2},$ &  & $x_{k}=x_{k-1}(1-x_{k-1})$ & for $k\geq1$%
\end{tabular}
\
\]
arose as an accident in \cite{Fi-logistic}, enabling us to understand the fine
details of $x_{k}\rightarrow0$ (as $k\rightarrow\infty$) more precisely than
ever imagined. Iterates of the logistic map occupy the memories of many
students, past and present. \ We examine here two cubic analogs of this -- not
as an attempt to generalize on a vast scale -- but simply as a means to
explore the material. \ For each example, two initial conditions are examined
(not just one as before). \ Also, in an interlude, another example is viewed
that is only loosely associated with the logistic map (if even that can be
said), but it is a cubic and it is well worth contemplating.

\section{The Map $x\mapsto x(1-x)^{2}$}

Given $0<x_{0}<1$ and%
\[%
\begin{array}
[c]{ccc}%
x_{k}=x_{k-1}\left(  1-x_{k-1}\right)  ^{2} &  & \text{for }k\geq1
\end{array}
\]
we start with asymptotics \cite{St-logistic, IS-logistic}%
\[
x_{k}\sim\frac{1}{2k}-\frac{3}{8}\frac{\ln(k)}{k^{2}}-\frac{C}{k^{2}}%
\]
valid as $k\rightarrow\infty$, for some constant $C=C(x_{0})$. \ On the basis
of numerical experimentation, we hypothesize that the next terms of the
asymptotic series must be of the form%
\[
u\,\frac{\ln(k)^{2}}{k^{3}}+v\,\frac{\ln(k)}{k^{3}}+\frac{w}{k^{3}}%
+p\,\frac{\ln(k)^{3}}{k^{4}}+q\,\frac{\ln(k)^{2}}{k^{4}}+r\,\frac{\ln
(k)}{k^{4}}+\frac{s}{k^{4}}.
\]
The challenge is to express each coefficient $u$, $v$, $w$, $p$, $q$, $r$, $s$
as a polynomial in $C$. \ To find these, we replace $k$ by $k+1$ everywhere:%
\[
\frac{1}{2(k+1)}-\frac{3}{8}\frac{\ln(k+1)}{(k+1)^{2}}-\frac{C}{(k+1)^{2}%
}+u\,\frac{\ln(k+1)^{2}}{(k+1)^{3}}+\cdots+\frac{s}{(k+1)^{4}}%
\]
and expand in powers of $k$ and $\ln(k)$:%
\[
\frac{1}{k+1}\sim\frac{1}{k}-\frac{1}{k^{2}}+\frac{1}{k^{3}}-\frac{1}{k^{4}%
}+\frac{1}{k^{5}}-+\cdots,
\]%
\[
\frac{\ln(k+1)}{(k+1)^{2}}\sim\left(  \frac{1}{k^{2}}-\frac{2}{k^{3}}+\frac
{3}{k^{4}}-\frac{4}{k^{5}}+-\cdots\right)  \ln(k)+\left(  \frac{1}{k^{3}%
}-\frac{5}{2k^{4}}+\frac{13}{3k^{5}}-+\cdots\right)  ,
\]%
\[
\frac{1}{(k+1)^{2}}\sim\frac{1}{k^{2}}-\frac{2}{k^{3}}+\frac{3}{k^{4}}%
-\frac{4}{k^{5}}+-\cdots,
\]%
\[
\frac{\ln(k+1)^{2}}{(k+1)^{3}}\sim\left(  \frac{1}{k^{3}}-\frac{3}{k^{4}%
}+\frac{6}{k^{5}}-+\cdots\right)  \ln(k)^{2}+\left(  \frac{2}{k^{4}}-\frac
{7}{k^{5}}+-\cdots\right)  \ln(k)+\left(  \frac{1}{k^{5}}-+\cdots\right)  ,
\]%
\[
\frac{\ln(k+1)}{(k+1)^{3}}\sim\left(  \frac{1}{k^{3}}-\frac{3}{k^{4}}+\frac
{6}{k^{5}}-+\cdots\right)  \ln(k)+\left(  \frac{1}{k^{4}}-\frac{7}{2k^{5}%
}+-\cdots\right)  ,
\]%
\[
\frac{1}{(k+1)^{3}}\sim\frac{1}{k^{3}}-\frac{3}{k^{4}}+\frac{6}{k^{5}}%
-+\cdots,
\]%
\[
\frac{\ln(k+1)^{3}}{(k+1)^{4}}\sim\left(  \frac{1}{k^{4}}-\frac{4}{k^{5}%
}+-\cdots\right)  \ln(k)^{3}+\left(  \frac{3}{k^{5}}-+\cdots\right)
\ln(k)^{2},
\]%
\[
\frac{\ln(k+1)^{2}}{(k+1)^{4}}\sim\left(  \frac{1}{k^{4}}-\frac{4}{k^{5}%
}+-\cdots\right)  \ln(k)^{2}+\left(  \frac{2}{k^{5}}-+\cdots\right)  \ln(k),
\]%
\[
\frac{\ln(k+1)}{(k+1)^{4}}\sim\left(  \frac{1}{k^{4}}-\frac{4}{k^{5}}%
+-\cdots\right)  \ln(k)+\left(  \frac{1}{k^{5}}-+\cdots\right)  ,
\]%
\[
\frac{1}{(k+1)^{4}}\sim\frac{1}{k^{4}}-\frac{4}{k^{5}}+-\cdots.
\]
Upon rearrangement, $x_{k+1}$ becomes%
\begin{align*}
&  -4p\frac{\ln(k)^{3}}{k^{5}}+\left(  6u+3p-4q\right)  \frac{\ln(k)^{2}%
}{k^{5}}+\left(  \frac{3}{2}-7u+6v+2q-4r\right)  \frac{\ln(k)}{k^{5}}\\
&  +\left(  -\frac{9}{8}+4C+u-\frac{7}{2}v+6w+r-4s\right)  \frac{1}{k^{5}%
}+p\,\frac{\ln(k)^{3}}{k^{4}}+\left(  -3u+q\right)  \frac{\ln(k)^{2}}{k^{4}}\\
&  +\left(  -\frac{9}{8}+2u-3v+r\right)  \frac{\ln(k)}{k^{4}}+\left(  \frac
{7}{16}-3C+v-3w+s\right)  \frac{1}{k^{4}}+u\,\frac{\ln(k)^{2}}{k^{3}}\\
&  +\left(  \frac{3}{4}+v\right)  \frac{\ln(k)}{k^{3}}+\left(  \frac{1}%
{8}+2C+w\right)  \frac{1}{k^{3}}-\frac{3}{8}\frac{\ln(k)}{k^{2}}-\left(
\frac{1}{2}+C\right)  \frac{1}{k^{2}}+\frac{1}{2k}.
\end{align*}
Performing an analogous substitution in $x_{k}$, the expression $x_{k}%
-2x_{k}^{2}+x_{k}^{3}$ becomes%
\begin{align*}
&  \left(  \frac{3}{2}u-2p\right)  \frac{\ln(k)^{3}}{k^{5}}+\left(  \frac
{27}{128}+\frac{3}{4}u+4C\,u+\frac{3}{2}v-2q\right)  \frac{\ln(k)^{2}}{k^{5}%
}\\
&  +\left(  \frac{9}{8}C+\frac{3}{4}v+4C\,v+\frac{3}{2}w-2r\right)  \frac
{\ln(k)}{k^{5}}+\left(  \frac{3}{2}C^{2}+\frac{3}{4}w+4C\,w-2s\right)
\frac{1}{k^{5}}\\
&  +p\,\frac{\ln(k)^{3}}{k^{4}}+\left(  -\frac{9}{32}-2u+q\right)  \frac
{\ln(k)^{2}}{k^{4}}+\left(  -\frac{9}{32}-\frac{3}{2}C-2v+r\right)  \frac
{\ln(k)}{k^{4}}\\
&  +\left(  -\frac{3}{4}C-2C^{2}-2w+s\right)  \frac{1}{k^{4}}+u\,\frac
{\ln(k)^{2}}{k^{3}}\\
&  +\left(  \frac{3}{4}+v\right)  \frac{\ln(k)}{k^{3}}+\left(  \frac{1}%
{8}+2C+w\right)  \frac{1}{k^{3}}-\frac{3}{8}\frac{\ln(k)}{k^{2}}-\left(
\frac{1}{2}+C\right)  \frac{1}{k^{2}}+\frac{1}{2k}.
\end{align*}
Matching coefficients, we obtain%

\[%
\begin{array}
[c]{ccccc}%
u=\dfrac{9}{32}, &  & v=\dfrac{3}{2}C-\dfrac{9}{32}, &  & w=2C^{2}-\dfrac
{3}{4}C+\dfrac{5}{32}%
\end{array}
\]
which are consistent with \cite{P1-logistic, P2-logistic} and%
\[%
\begin{array}
[c]{ccccc}%
p=-\dfrac{27}{128}, &  & q=-\left(  \dfrac{27}{16}C-\dfrac{135}{256}\right)
, &  & r=-\left(  \dfrac{9}{2}C^{2}-\dfrac{45}{16}C+\dfrac{9}{16}\right)  ,
\end{array}
\]%
\[
s=-\left(  4C^{3}-\dfrac{15}{4}C^{2}+\dfrac{3}{2}C-\dfrac{51}{256}\right)
\]
which are new (as far as is known).

These seven parameter values allow us to estimate the constant $C$. \ We do
this both for $x_{0}=1/2$ (the midpoint) and for $x_{0}=1/3$ (the argument at
which $C$ is minimal). \ Our simple procedure involves computing
$a_{100000000}$ exactly via recursion, setting this equal to our series (up to
$s/k^{4}$) and then solving:%
\[%
\begin{array}
[c]{ccc}%
C\left(  \frac{1}{2}\right)  =1.025737030693254..., &  & C\left(  \frac{1}%
{3}\right)  =0.787336122933677....
\end{array}
\]
Note that the estimate $8\,C(1/2)\approx8.205896246$ appears in
\cite{S1-logistic}. \ Also, a simple change of variable would permit one to
analyze $x_{k}=x_{k-1}^{2}\left(  1-x_{k-1}\right)  $.

\section{An Interlude}

The cubic recurrence%
\[%
\begin{tabular}
[c]{llll}%
$x_{0}=0,$ &  & $x_{k}=\dfrac{2}{3}+\dfrac{1}{3}\,x_{k-1}^{3}$ & for $k\geq1$%
\end{tabular}
\ \ \
\]
arises in the study of random Galton-Watson ternary\ tree heights at
criticality. \ Its binary analog \cite{Fi-logistic}%
\[%
\begin{tabular}
[c]{llll}%
$\xi_{0}=0,$ &  & $\xi_{k}=\dfrac{1}{2}+\dfrac{1}{2}\,\xi_{k-1}^{2}$ & for
$k\geq1$%
\end{tabular}
\ \ \
\]
is easily transformed to the logistic%
\[%
\begin{tabular}
[c]{llll}%
$\eta_{0}=\dfrac{1}{2},$ &  & $\eta_{k}=\eta_{k-1}(1-\eta_{k-1})$ & for
$k\geq1$%
\end{tabular}
\ \ \
\]
under $\xi_{k}=1-2\eta_{k}$. \ Letting instead $x_{k}=1-3y_{k}$, we obtain%
\begin{align*}
1-3y_{k}  &  =\frac{1}{3}\left[  2+\left(  1-3y_{k-1}\right)  ^{3}\right] \\
&  =\frac{1}{3}\left(  3-9y_{k-1}+27y_{k-1}^{2}-27y_{k-1}^{3}\right)
=1-3y_{k-1}+9y_{k-1}^{2}-9y_{k-1}^{3}%
\end{align*}
hence%
\[
y_{k}=y_{k-1}-3y_{k-1}^{2}+3y_{k-1}^{3}=y_{k-1}\left(  1-3y_{k-1}+3y_{k-1}%
^{2}\right)  .
\]
This expression does not resemble a logistic; it nevertheless allows starting
asymptotics \cite{St-logistic, IS-logistic}%
\[
\frac{1-x_{k}}{3}=y_{k}\sim\frac{1}{3k}-\frac{2}{9}\frac{\ln(k)}{k^{2}}%
-\frac{C}{3k^{2}}%
\]
valid as $k\rightarrow\infty$, for some constant $C$. \ Thus%
\[
x_{k}\sim1-\frac{1}{k}+\frac{2}{3}\frac{\ln(k)}{k^{2}}+\frac{C}{k^{2}}.
\]
On the basis of numerical experimentation, we hypothesize that the next terms
of the asymptotic series must be of the form%
\[
u\,\frac{\ln(k)^{2}}{k^{3}}+v\,\frac{\ln(k)}{k^{3}}+\frac{w}{k^{3}}%
+p\,\frac{\ln(k)^{3}}{k^{4}}+q\,\frac{\ln(k)^{2}}{k^{4}}+r\,\frac{\ln
(k)}{k^{4}}+\frac{s}{k^{4}}.
\]
The algorithm for finding coefficients $u$, $v$, $w$, $p$, $q$, $r$, $s$ is as
before. \ Upon rearrangement, the terms of $x_{k+1}$ involving either $k^{4}$
or $k^{5}$ become most relevant. \ Performing an analogous substitution in
$x_{k}$, we match coefficients from the expression $2/3+x_{k-1}^{3}/3$ with
those of $x_{k+1}$. \ The terms containing $k^{4}$ give%
\[%
\begin{array}
[c]{ccccc}%
u=-\dfrac{4}{9}, &  & v=-\left(  \dfrac{4}{3}C-\dfrac{4}{9}\right)  , &  &
w=-\left(  C^{2}-\dfrac{2}{3}C+\dfrac{2}{9}\right)
\end{array}
\]
which are consistent with \cite{P1-logistic, P2-logistic}. \ The terms
containing $k^{5}$ give%
\[%
\begin{array}
[c]{cccc}%
p=\dfrac{1}{3}, & q=\dfrac{4}{3}C-\dfrac{37}{54}, &  & r=2C^{2}-\dfrac{20}%
{9}C+\dfrac{43}{54},
\end{array}
\]%
\[
s=C^{3}-\frac{5}{3}C^{2}+\frac{10}{9}C-\frac{23}{108}%
\]
which are new (as far as is known).

These seven parameter values allow us to estimate the constant $C$. \ We do
this only for $x_{0}=0$. \ Our simple procedure involves computing
$a_{100000000}$ exactly via recursion, setting this equal to our series (up to
$s/k^{4}$) and then solving:%
\[
C=1.354567323982625....
\]
Similar analysis of $x_{k}=3/4+x_{k-1}^{4}/4$ is possible (for
random\ Galton-Watson quaternary tree heights at criticality).

\section{The Map $x\mapsto x(1-x^{2})$}

Given $0<x_{0}<1$ and%
\[%
\begin{array}
[c]{ccc}%
x_{k}=x_{k-1}\left(  1-x_{k-1}^{2}\right)  &  & \text{for }k\geq1
\end{array}
\]
we start with asymptotics \cite{St-logistic, IS-logistic}%
\[
x_{k}\sim\frac{1}{\sqrt{2}}\frac{1}{k^{1/2}}-\frac{3}{8\sqrt{2}}\frac{\ln
(k)}{k^{3/2}}-\frac{C}{2\sqrt{2}}\frac{1}{k^{3/2}}%
\]
valid as $k\rightarrow\infty$, for some constant $C$. \ It is reasonable to
hypothesize that the next terms of this series are%
\[
u\,\frac{\ln(k)^{2}}{k^{5/2}}+v\,\frac{\ln(k)}{k^{5/2}}+\frac{w}{k^{5/2}%
}+p\,\frac{\ln(k)^{3}}{k^{7/2}}+q\,\frac{\ln(k)^{2}}{k^{7/2}}+r\,\frac{\ln
(k)}{k^{7/2}}+\frac{s}{k^{7/2}}.
\]
To find $u$, $v$, $w$, $p$, $q$, $r$, $s$, we replace $k$ by $k+1$
everywhere:
\[
\frac{1}{\sqrt{2}}\frac{1}{(k+1)^{1/2}}-\frac{3}{8\sqrt{2}}\frac{\ln
(k+1)}{(k+1)^{3/2}}-\frac{C}{2\sqrt{2}}\frac{1}{(k+1)^{3/2}}+u\,\frac
{\ln(k+1)^{2}}{(k+1)^{5/2}}+\cdots+\frac{s}{(k+1)^{7/2}}%
\]
and expand in powers of $k$ and $\ln(k)$:%
\[
\frac{1}{(k+1)^{1/2}}\sim\frac{1}{k^{1/2}}-\frac{1}{2k^{3/2}}+\frac
{3}{8k^{5/2}}-\frac{5}{16k^{7/2}}+\frac{35}{128k^{9/2}}-+\cdots,
\]%
\[
\frac{\ln(k+1)}{(k+1)^{3/2}}\sim\left(  \frac{1}{k^{3/2}}-\frac{3}{2k^{5/2}%
}+\frac{15}{8k^{7/2}}-\frac{35}{16k^{9/2}}+-\cdots\right)  \ln(k)+\left(
\frac{1}{k^{5/2}}-\frac{2}{k^{7/2}}+\frac{71}{24k^{9/2}}-+\cdots\right)  ,
\]%
\[
\frac{1}{(k+1)^{3/2}}\sim\frac{1}{k^{3/2}}-\frac{3}{2k^{5/2}}+\frac
{15}{8k^{7/2}}-\frac{35}{16k^{9/2}}+-\cdots,
\]%
\[
\frac{\ln(k+1)^{2}}{(k+1)^{5/2}}\sim\left(  \frac{1}{k^{5/2}}-\frac
{5}{2k^{7/2}}+\frac{35}{8k^{9/2}}-+\cdots\right)  \ln(k)^{2}+\left(  \frac
{2}{k^{7/2}}-\frac{6}{k^{9/2}}+-\cdots\right)  \ln(k)+\left(  \frac{1}%
{k^{9/2}}-+\cdots\right)  ,
\]%
\[
\frac{\ln(k+1)}{(k+1)^{5/2}}\sim\left(  \frac{1}{k^{5/2}}-\frac{5}{2k^{7/2}%
}+\frac{35}{8k^{9/2}}-+\cdots\right)  \ln(k)+\left(  \frac{1}{k^{7/2}}%
-\frac{3}{k^{9/2}}+-\cdots\right)  ,
\]%
\[
\frac{1}{(k+1)^{5/2}}\sim\frac{1}{k^{5/2}}-\frac{5}{2k^{7/2}}+\frac
{35}{8k^{9/2}}-+\cdots,
\]%
\[
\frac{\ln(k+1)^{3}}{(k+1)^{7/2}}\sim\left(  \frac{1}{k^{7/2}}-\frac
{7}{2k^{9/2}}+-\cdots\right)  \ln(k)^{3}+\left(  \frac{3}{k^{9/2}}%
-+\cdots\right)  \ln(k)^{2},
\]%
\[
\frac{\ln(k+1)^{2}}{(k+1)^{7/2}}\sim\left(  \frac{1}{k^{7/2}}-\frac
{7}{2k^{9/2}}+-\cdots\right)  \ln(k)^{2}+\left(  \frac{2}{k^{9/2}}%
-+\cdots\right)  \ln(k),
\]%
\[
\frac{\ln(k+1)}{(k+1)^{7/2}}\sim\left(  \frac{1}{k^{7/2}}-\frac{7}{2k^{9/2}%
}+-\cdots\right)  \ln(k)+\left(  \frac{1}{k^{9/2}}-+\cdots\right)  ,
\]%
\[
\frac{1}{(k+1)^{7/2}}\sim\frac{1}{k^{7/2}}-\frac{7}{2k^{9/2}}+-\cdots.
\]
Upon rearrangement, the terms of $x_{k+1}$ involving either $k^{7/2}$ or
$k^{9/2}$ become most relevant. \ Performing an analogous substitution in
$x_{k}$, we match coefficients from the expression $x_{k}-x_{k}^{3}$ with
those of $x_{k+1}$. \ The terms containing $k^{7/2}$ give%
\[%
\begin{array}
[c]{ccccc}%
u=\dfrac{27}{128\sqrt{2}}, &  & v=\dfrac{\frac{9}{16}C-\frac{9}{32}}{\sqrt{2}%
}, &  & w=\dfrac{\frac{3}{8}C^{2}-\frac{3}{8}C+\frac{5}{32}}{\sqrt{2}}%
\end{array}
\]
which are consistent with \cite{LSZ-logistic}. \ The terms containing
$k^{9/2}$ give%
\[%
\begin{array}
[c]{ccccc}%
p=-\dfrac{135}{1024\sqrt{2}}, &  & q=-\dfrac{\frac{135}{256}C-\frac{27}{64}%
}{\sqrt{2}}, &  & r=-\dfrac{\frac{45}{64}C^{2}-\frac{9}{8}C+\frac{129}{256}%
}{\sqrt{2}},
\end{array}
\]%
\[
s=-\dfrac{\frac{5}{16}C^{3}-\frac{3}{4}C^{2}+\frac{43}{64}C-\frac{51}{256}%
}{\sqrt{2}}%
\]
which are new (as far as is known).

These seven parameter values allow us to estimate the constant $C$. \ We do
this both for $x_{0}=1/2$ (the midpoint) and for $x_{0}=1/\sqrt{3}$ (the
argument at which $C$ is minimal). \ Our simple procedure involves computing
$a_{100000000}$ exactly via recursion, setting this equal to our series (up to
$s/k^{7/2}$) and then solving:%
\[%
\begin{array}
[c]{ccc}%
C\left(  \frac{1}{2}\right)  =1.709211474227594..., &  & C\left(  \frac
{1}{\sqrt{3}}\right)  =1.574672245867354....
\end{array}
\]
The latter improves upon an estimate $1.5739$ reported in \cite{LSZ-logistic}.

\section{An Impromptu}

Letting%
\[
\eta_{k}=\frac{1}{a_{k}}=\frac{1}{1+b_{k}}%
\]
where $\eta_{k}$ satisfies the logistic recurrence (see the interlude), we
have \cite{IS-logistic, GM0-logistic}%
\[%
\begin{array}
[c]{ccc}%
a_{k+1}=a_{k}+1+\dfrac{1}{a_{k}-1}, &  & b_{k+1}=1+b_{k}+\dfrac{1}{b_{k}}%
\end{array}
\]
because $b_{k}=a_{k}-1=1/\eta_{k}-1$,%
\[
\frac{1}{\eta_{k}}+1+\frac{1}{\frac{1}{\eta_{k}}-1}=\frac{(1-\eta_{k}%
)+\eta_{k}(1-\eta_{k})+\eta_{k}^{2}}{\eta_{k}(1-\eta_{k})}=\frac{1}{\eta
_{k+1}},
\]%
\[
1+\frac{1-\eta_{k}}{\eta_{k}}+\frac{\eta_{k}}{1-\eta_{k}}=\frac{\eta
_{k}(1-\eta_{k})+(1-\eta_{k})^{2}+\eta_{k}^{2}}{\eta_{k}(1-\eta_{k})}%
=\frac{1-\eta_{k}(1-\eta_{k})}{\eta_{k}(1-\eta_{k})}=\frac{1-\eta_{k+1}}%
{\eta_{k+1}}.
\]
Examining these formulas side-by-side, it is natural to inquire about%
\[
x_{k+1}=x_{k}+\dfrac{1}{x_{k}}%
\]
where $x_{0}>0$. \ Here the sequence grows without bound and its asymptotic
expansion begins as \cite{P1-logistic, P2-logistic}%
\begin{align*}
&  \sqrt{2k}+\frac{1}{4\sqrt{2}}\frac{\ln(k)}{k^{1/2}}+\frac{C}{k^{1/2}}%
-\frac{1}{64\sqrt{2}}\frac{\ln(k)^{2}}{k^{3/2}}-\frac{4C-\sqrt{2}}{32}%
\,\frac{\ln(k)}{k^{3/2}}-\frac{8\sqrt{2}C^{2}-8C+\sqrt{2}}{32}\,\frac
{1}{k^{3/2}}\\
&  +\frac{1}{512\sqrt{2}}\frac{\ln(k)^{3}}{k^{5/2}}+\frac{3C-\sqrt{2}}%
{128}\,\frac{\ln(k)^{2}}{k^{5/2}}+\frac{24\sqrt{2}C^{2}-32C+5\sqrt{2}}%
{256}\,\frac{\ln(k)}{k^{5/2}}\\
&  +\frac{192C^{3}-192\sqrt{2}C^{2}+120C-11\sqrt{2}}{768}\,\frac{1}{k^{5/2}}.
\end{align*}
With initial condition $x_{0}=1$, we determine the constant
$C=0.609222829204782...$, omitting details. \ Note that the estimate
$C/\sqrt{2}\approx0.43078559$ appears in \cite{S2-logistic}; two similar
iterative schemes appear in \cite{S3-logistic, S4-logistic}.

At first glance, one might confuse the above with a widely known recurrence
for continued fractions:%
\[
x_{k+1}=x_{0}+\dfrac{1}{x_{k}}%
\]
which is far simpler. \ When $x_{0}=2$, it yields%
\[
\lim_{k\rightarrow\infty}\,x_{k}=1+\sqrt{2}=\frac{1}{-1+\sqrt{2}}%
\]
and the convergence is exponential. \ A generalization of this involves what
are called bifurcating continued fractions \cite{Fn-logistic, GM1-logistic,
GM2-logistic}:%
\[
y_{k+1}=y_{0}+\dfrac{y_{0}+\tfrac{1}{y_{k}}}{y_{k}}%
\]
and, when $y_{0}=3$, it yields%
\[
\lim_{k\rightarrow\infty}\,y_{k}=1+2^{1/3}+2^{2/3}=\frac{1}{-1+2^{1/3}}.
\]
This little diversion might serve as inspiration for someone to examine
\cite{P1-logistic, P2-logistic}%
\[
z_{k+1}=z_{k}+\dfrac{1}{z_{k}^{2}}%
\]
where $z_{0}>0$. \ Demonstrably $z_{k}\sim(3k)^{1/3}$, followed by a term of
order $\ln(k)/k^{2/3}$ and after this (once again) a term $C/k^{2/3}$.

In summary, high-precision calculation of the constant $C$, given an
associated nonlinear recurrence $x_{k+1}=f(x_{k})$ and a starting value
$x_{0}$, is challenging. \ There are trivial cases that needn't be considered,
e.g.,
\[
x_{k+1}=x_{k}^{2}%
\]
(exponential growth or decay) for which%
\[
x_{k}=x_{0}^{2^{k}}%
\]
and%
\[
x_{k+1}=\frac{1}{x_{0}+\tfrac{1}{x_{k}}}=\frac{x_{k}}{1+x_{0}x_{k}}%
\]
(reciprocal of continued fractions) for which%
\[
x_{k}=\frac{x_{0}}{1+k\,x_{0}^{2}}.
\]
Among nontrivial cases, we wonder whether there exists an explicit recurrence
$f$ and an explicit value $x_{0}$ such that the corresponding constant $C$ is
known exactly (i.e., possesses a closed-form expression). \ A resolution of
this issue would be good to learn someday.

\section{A\ Reprise}

The sequence%
\[%
\begin{array}
[c]{ccccc}%
x_{k+1}=x_{k}+\dfrac{a}{x_{k}}, &  & a>0, &  & x_{0}=b>0
\end{array}
\]
generalizes our earlier work (on the special case $a=b=1$) and its asymptotic
expansion begins as%
\begin{align*}
&  \sqrt{2a\,k}+\frac{\sqrt{a}}{4\sqrt{2}}\frac{\ln(k)}{k^{1/2}}+\frac
{C}{k^{1/2}}-\frac{\sqrt{a}}{64\sqrt{2}}\frac{\ln(k)^{2}}{k^{3/2}}%
-\frac{4C-\sqrt{2a}}{32}\,\frac{\ln(k)}{k^{3/2}}-\frac{8\sqrt{2}C^{2}%
-8C+\sqrt{2a}}{32}\,\frac{1}{k^{3/2}}\\
&  +\frac{\sqrt{a}}{512\sqrt{2}}\frac{\ln(k)^{3}}{k^{5/2}}+\frac{3C-\sqrt{2a}%
}{128}\,\frac{\ln(k)^{2}}{k^{5/2}}+\frac{4\sqrt{2}\left(  \sqrt{a}+5\right)
C^{2}-32\sqrt{a}C+5\sqrt{2}a}{256\sqrt{a}}\,\frac{\ln(k)}{k^{5/2}}\\
&  +\frac{96\left(  \sqrt{a}+1\right)  C^{3}-2\sqrt{2a}\left(  87\sqrt
{a}+9\right)  C^{2}+120a\,C-11\sqrt{2}a^{3/2}}{768a}\,\frac{1}{k^{5/2}}.
\end{align*}
We tabulate values of $C=C(a,b)$ for $a=\frac{1}{2},1,2$ and $b=\frac{1}%
{\sqrt{2}},1,\sqrt{2}$. \ For fixed $a$,\ $\min_{b>0}C(a,b)$ is starred.
$\medskip$\ 

\textit{Table of }$C(a,b)$\textit{ values for selected }$a$ \textit{and} $b$
\
\[%
\begin{tabular}
[c]{|l|l|l|l|}\hline
$a\backslash b$ & $1/\sqrt{2}$ & $1$ & $\sqrt{2}$\\\hline
$1/2$ & $0.430785593784355...^{\ast}$ & $0.538934324436812...$ &
$0.930785593784355...$\\\hline
$1$ & $0.762168230846922...$ & $0.609222829204782...^{\ast}$ &
$0.762168230846922...$\\\hline
$2$ & $1.861571187568711...$ & $1.077868648873625...$ &
$0.861571187568711...^{\ast}$\\\hline
\end{tabular}
\ \ \ \ \ \ \ \
\]
Elements along the diagonals (northwest to southeast) can be related via%
\[%
\begin{array}
[c]{ccc}%
C\left(  a,\sqrt{a}\right)  =\sqrt{a}\,C(1,1), &  & C(1,b)=b\,C\left(
\dfrac{1}{b^{2}},1\right)
\end{array}
\]
and diametrically opposed elements via
\[
C(a,b)=a\,C\left(  \dfrac{1}{a},\frac{1}{b}\right)  .
\]
Let's prove the latter:\ if $u_{0}=b$ and $v_{0}=\frac{1}{b}$, then%
\[%
\begin{array}
[c]{ccccc}%
u_{1}=b+\dfrac{a}{b} &  & \text{and} &  & v_{1}=\dfrac{1}{b}+\dfrac{\frac
{1}{a}}{\frac{1}{b}}=\dfrac{1}{b}+\dfrac{b}{a}=\dfrac{1}{a}\left(  b+\dfrac
{a}{b}\right)  =\dfrac{u_{1}}{a}.
\end{array}
\]
Assuming $v_{k-1}=\frac{u_{k-1}}{a}$ for some $k\geq2$, we have%
\[%
\begin{array}
[c]{ccccc}%
u_{k}=u_{k-1}+\dfrac{a}{u_{k-1}} &  & \text{and} &  & v_{k}=v_{k-1}%
+\dfrac{\frac{1}{a}}{v_{k-1}}=\dfrac{u_{k-1}}{a}+\dfrac{1}{u_{k-1}}=\dfrac
{1}{a}\left(  u_{k-1}+\dfrac{a}{u_{k-1}}\right)  =\dfrac{u_{k}}{a}%
\end{array}
\]
and the formula is true by induction. \ Let's briefly examine the former:\ if
$u_{0}=\sqrt{a}$ and $v_{0}=1$, then%
\[%
\begin{array}
[c]{ccccc}%
u_{1}=\sqrt{a}+\dfrac{a}{\sqrt{a}}=2\sqrt{a} &  & \text{and} &  &
v_{1}=1+\dfrac{1}{1}=2
\end{array}
\]
and if $u_{0}=b$ and $v_{0}=1$, then%
\[%
\begin{array}
[c]{ccccc}%
u_{1}=b+\dfrac{1}{b}=b\left(  1+\dfrac{1}{b^{2}}\right)  &  & \text{and} &  &
v_{1}=1+\dfrac{\frac{1}{b^{2}}}{1}=1+\dfrac{1}{b^{2}}.
\end{array}
\]

Elements on the top and bottom rows satisfy%
\[%
\begin{array}
[c]{ccc}%
C\left(  \dfrac{1}{2},\dfrac{1}{\sqrt{2}}\right)  =C\left(  \dfrac{1}{2}%
,\sqrt{2}\right)  -\dfrac{1}{2}, &  & C\left(  2,\dfrac{1}{\sqrt{2}}\right)
=C\left(  2,\sqrt{2}\right)  +1;
\end{array}
\]
we'll prove the latter by noticing that $u_{k}=v_{k+1}$ for all $k\geq1$ and%
\[
u_{k}\sim2\sqrt{k}+\frac{1}{4}\frac{\ln(k)}{\sqrt{k}}+\frac{C_{u}}{\sqrt{k}},
\]%
\[
v_{k+1}\sim2\sqrt{k+1}+\frac{1}{4}\frac{\ln(k+1)}{\sqrt{k+1}}+\frac{C_{v}%
}{\sqrt{k+1}}\sim2\sqrt{k}+\frac{1-\frac{1}{4k}}{\sqrt{k}}+\frac{1}{4}%
\frac{\left(  1-\frac{1}{2k}\right)  \ln(k)}{\sqrt{k}}+\frac{\left(
1-\frac{1}{2k}\right)  C_{v}}{\sqrt{k}}%
\]
thus $0=C_{u}-1-C_{v}$ upon taking differences. \ Let's likewise examine the
former:\ here $u_{k+1}=v_{k}$ for all $k\geq1$, thus%
\[
v_{k}\sim\sqrt{k}+\frac{1}{8}\frac{\ln(k)}{\sqrt{k}}+\frac{C_{v}}{\sqrt{k}},
\]%
\[
u_{k+1}\sim\sqrt{k+1}+\frac{1}{8}\frac{\ln(k+1)}{\sqrt{k+1}}+\frac{C_{u}%
}{\sqrt{k+1}}\sim\sqrt{k}+\frac{1-\frac{1}{4k}}{2\sqrt{k}}+\frac{1}{8}%
\frac{\left(  1-\frac{1}{2k}\right)  \ln(k)}{\sqrt{k}}+\frac{\left(
1-\frac{1}{2k}\right)  C_{u}}{\sqrt{k}}%
\]
and $0=C_{v}-\frac{1}{2}-C_{u}$.

In summary, the table can be completely generated by the two (algebraically
independent?) constants $C(1,1)$ and $C\left(  1,1/\sqrt{2}\right)  $. \ 

Note that the estimate $\frac{1}{2}C(2,1)\approx0.538934322$ appears in
\cite{S3-logistic} and $2C\left(  \frac{1}{2},\sqrt{2}\right)  -\frac{1}%
{2}\approx1.36157$ appears in \cite{S4-logistic}. \ A\ new reference
(discovered after Section 4 was written) involves $C(1,1)$ and gives two
numerical evaluations to high precision \cite{S5-logistic, S6-logistic}:%
\[
2\sqrt{2}\,C-2=-0.2768576248625765389364372...,
\]%
\[
\exp\left(  4\sqrt{2}\,C-4\right)  =0.5748102746717850652963478...
\]
from which we deduce%
\[
C(1,1)=0.6092228292047829402293060....
\]
A\ similarly accurate evaluation of $C\left(  1,1/\sqrt{2}\right)  $ would be
most welcome.

As an aside, Newton's method for solving the equation $z^{2}=c$ gives
recurrences \cite{GS-logistic}%
\[%
\begin{array}
[c]{ccccccc}%
x_{n+1}=\dfrac{1}{2}x_{n}+\dfrac{1}{x_{n}}, &  & x_{0}=1; &  & y_{n+1}%
=\dfrac{1}{2}y_{n}+\dfrac{3}{2}\dfrac{1}{y_{n}}, &  & y_{0}=2
\end{array}
\]
for $c=2$ \&\ $3$ and asymptotics \cite{S7-logistic, S8-logistic}%
\[%
\begin{array}
[c]{ccc}%
\left(  x_{n}-\sqrt{2}\right)  \left(  3+2\sqrt{2}\right)  ^{2^{n}}%
\rightarrow2\sqrt{2}; &  & \left(  y_{n}-\sqrt{3}\right)  \left(  7+4\sqrt
{3}\right)  ^{2^{n}}\rightarrow2\sqrt{3}%
\end{array}
\]
as $n\rightarrow\infty$. \ Solving instead $1/z^{2}=c$ gives division-free
recurrences \cite{BC-logistic}%
\[%
\begin{array}
[c]{ccccccc}%
x_{n+1}=x_{n}\left(  \dfrac{3}{2}-x_{n}^{2}\right)  , &  & x_{0}=\dfrac{1}%
{2}; &  & y_{n+1}=\dfrac{3}{2}y_{n}\left(  1-y_{n}^{2}\right)  , &  &
y_{0}=\dfrac{1}{3}%
\end{array}
\]
but exact asymptotics (akin to the preceding) remain open \cite{S9-logistic,
S0-logistic}.

\section{Acknowledgements}

I\ am thankful to Yong-Guo Shi for helpful discussion regarding
\cite{LSZ-logistic} and their newly updated estimate $1.57468$ of
$C(1/\sqrt{3})$. \ My attempts to reach Dumitru Popa \cite{P1-logistic,
P2-logistic} have regrettably failed. \ In the statement of Theorem 10 on page
24 of \cite{P2-logistic}, the lead coefficient $2a_{1}^{4}K^{2}$ of the
$1/n^{3}$ term should be $2a_{1}^{6}K^{2}$ (Popa's corresponding expression
$2a_{1}^{2}C^{2}$ on page 25, where $C=-a_{1}^{2}K$, is correct). \ This
algebraic error affects his Corollary 11. \ I\ deeply appreciate as well my
correspondence with Eugen Ionascu, Pante Stanica and Michael Somos.

\end{document}